\documentclass[final]{siamltex}
\usepackage{amssymb}
\usepackage{amsmath}
\usepackage{color}
\newtheorem{remark}[theorem]{Remark}
\newtheorem{claim}[theorem]{Claim}
\newcommand{\R}{\mbox{${\rm I}\!{\rm R}$}}

\newcommand{\D}{\Delta}
\newcommand{\I}{\mathcal I}
\newcommand{\J}{\mathcal J}
\newcommand{\DD}{{\D_m\times \D_n}}

\newcommand{\Oh}{{\mathcal{O}}}
\newcommand{\e}{\varepsilon}
\renewcommand{\k}{\kappa\,}
\renewcommand{\phi}{\varphi}
\def\Limsup{\mathop{{\rm Lim\,sup\,}}}
\def\co{\mathop{{\rm co\,}}}
\def\cone{\mathop{{\rm cone\,}}}
\newcommand{\dist}{{\rm \,dist\,}}
\newcommand{\dom}{{\rm dom\,}}
\newcommand{\epi}{{\rm epi\,}}
\newcommand{\lip}{{\rm lip\,}}
\newcommand{\reg}{{\rm reg\,}}

\newcommand{\rge}{{\rm rge\,}}
\newcommand{\gph}{{\rm gph\,}}

\newcommand{\cl}{{\rm cl\,}}

\renewcommand{\span}{{\rm span\, }}
\newcommand{\matr}[1]{\begin{bmatrix} #1 \end{bmatrix}}

\newcommand{\dc}{{\partial}}

\def\transp{^{\rm T}}

\newcommand{\B}{{I\!\!B}}
\newcommand{\N}{{I\!\!N}}

\def\disp{\displaystyle}

\def\Limsup{\mathop{{\rm Lim}\,{\rm sup}}}

\def\tto{\;{\lower 1pt \hbox{$\rightarrow$}}\kern -10pt
\hbox{\raise 2pt \hbox{$\rightarrow$}}\;}
\def\Hat{\widehat}
\def\Tilde{\widetilde}
\def\Bar{\overline}
\def\ra{\rangle}
\def\la{\langle}
\def\ve{\varepsilon}
\def\B{I\!\!B}
\def\h{\hfill\Box}
\def\R{I\!\!R}
\def\N{I\!\!N}

\def\ox{\bar{x}}
\def\oy{\bar{y}}
\def\oz{\bar{z}}

\def\ow{\Bar{w}}

\def\k{\kappa}
\def\co{\mbox{\rm co}\,}
\def\cone{\mbox{\rm cone}\,}

\def\gph{\mbox{\rm gph}\,}
\def\epi{\mbox{\rm epi}\,}

\def\span{\mbox{\rm span}\,}

\def\dom{\mbox{\rm dom}\,}
\def\ker{\mbox{\rm ker}\,}

\def\lip{\mbox{\rm lip}\,}
\def\reg{\mbox{\rm reg}\,}

\def\cl*co{\mbox{\rm cl}^*\mbox{\rm co}\,}
\def\cl{\mbox{\rm cl}\,}

\def\h{\hfill\triangle}

\def\O{\Omega}

\def\ph{\varphi}

\def\emp{\emptyset}
\def\st{\stackrel}
\def\oR{\Bar{\R}}
\def\lm{\lambda}
\def\gg{\gamma}
\def\dd{\delta}

\def\tilde{\Tilde}

\def\N{I\!\!N}
\def \hs7{\hspace*{7pt}}

\renewcommand{\theequation}{\thesection.\arabic{equation}}
\begin{document}
\begin{center}
{\bf APPLYING METRIC REGULARITY TO COMPUTE A CONDITION MEASURE OF
A SMOOTHING ALGORITHM FOR MATRIX GAMES}\\[3ex]

BORIS S. MORDUKHOVICH\footnote{Wayne State University, Detroit,
MI, USA (boris@math.wayne.edu). Research of this author was partly
supported by the US National Science Foundation under grants
DMS-0603846 and DMS-1007132 and by the Australian Research Council
under grant DP-12092508.}, JAVIER F. PE\~NA\footnote{Carnegie Mellon
University, Pittsburgh, PA, USA (jfp@andrew.cmu.edu). Research of
this author was partly supported by the US National Science
Foundation under grant CCF-0830533.} and VERA
ROSHCHINA\footnote{Universidade de \'Evora, \'Evora, Portugal
(veraroshchina@gmail.com).}
\end{center}

\small {\bf Abstract.} We develop an approach of variational
analysis and generalized differentiation to conditioning issues
for two-person zero-sum matrix games. Our major results establish
precise relationships between a certain condition measure of the
smoothing first-order algorithm proposed by Gilpin et al. [{\em Proceedings
of the 23rd AAAI Conference (2008) pp. 75--82}]
and the exact bound of metric regularity for an associated
set-valued mapping. In this way we compute the aforementioned
condition measure in terms of the initial matrix game data.\\



{\bf Abbreviated title.} Condition measure for matrix games\\

\normalsize
\section{Introduction and formulation of main results}\label{sec:Intro}
\setcounter{equation}{0}

This paper is devoted to applications of advanced techniques in
variational analysis and generalized differentiation to the study
of conditioning in optimization. Our specific goal is to apply the
key notions and generalized differential characterizations of
Lipschitzian stability and metric regularity, fundamental in
variational analysis, to computing a certain condition measure of
the first-order smoothing algorithm proposed in
\cite{GilPenaSandh} to find approximate Nash equilibria of
two-person zero-sum matrix games.

To the best of our knowledge, applications of Lipschitzian
stability and metric regularity to numerical aspects of
optimization were initiated by Robinson in the 1970s; see, e.g.,
\cite{rob} and the references therein. In the complexity theory,
Renegar \cite{Renegar95a,Renegar95b} established relationships
between the rate of convergence of interior-point methods for
linear and conic convex programs and their ``distance to
ill-posedness" and related condition numbers. We refer the reader
to \cite{ADG,ABRS,DonLewRock,KK,LewisLukeMalick2009,PeRe} and
their bibliographies for more recent results in this direction for
various algorithms in convex and nonconvex optimization
problems.\vspace*{0.05in}

In \cite{GilPenaSandh}, a new condition measure was introduced to
evaluate the complexity of a first-order algorithm for solving a
two-person zero-sum game
\begin{equation}\label{eq:equilgeneral}
\min_{x\in Q_1}\max_{y\in Q_2} x\transp A y =
\max_{y\in Q_2} \min_{x\in Q_1}x\transp A y,
\end{equation}
where $A\in\R^{m\times n}$, where the symbol $\transp$ stands for
transposition, and where each of the sets $Q_1$ and $Q_2$ is
either a simplex (in the matrix game formulation) or a more
elaborate polytope (in the case of sequential games). Problems of
this type arise in many interesting applications; see, e.g.,
\cite{OsboR94,Roma62,Sten96,s} and the references therein.

It was shown in \cite{GilPenaSandh} that an iterative version of
Nesterov's first-order smoothing algorithm
\cite{Nesterov05,NesterovMP05} computes an $\e$-equilibrium point
(in the sense of Nash) for problem (\ref{eq:equilgeneral}) in
$\Oh(\|A\|\k(A)\ln (1/\e))$ iterations, where $\kappa(A)$ is a
condition measure of (\ref{eq:equilgeneral}) depending only on $A$;
see \eqref{eq:k} for the definition of the condition measure
$\k(A)$ in the case of matrix games. The dependence of this
complexity bound on $\e$ is exponentially better than the
complexity bound $\Oh(1/\e)$ in the original Nesterov's smoothing
techniques. Furthermore, it was proved in \cite{GilPenaSandh} that
the condition measure $\kappa(A)$ is always finite while the proof
therein was non-constructive. In particular, no explicit upper
bound on $\k(A)$ was given. However, numerical results reported
in~\cite{GilPenaSandh} clearly demonstrate that the developed
iterative version of Nesterov's smoothing algorithm is faster than
other algorithms known in this setting and that this version
exhibits at least linear convergence for a random collection of
the problem instances considered in \cite{GilPenaSandh}. This
allows us to treat the number $\kappa(A)$ as a condition measure
of the algorithm and evaluate it in what follows. \vspace*{0.05in}

In this paper we focus on the {\em matrix game equilibrium
problem}
\begin{equation}\label{eq:equil}
\min_{x\in\D_m}\max_{y\in \D_n}x\transp A y =\max_{y\in\D_n}
\min_{x\in\D_m}x\transp A y,
\end{equation}
where the $m$-dimensional simplex
\begin{eqnarray*}
\D_m:=\disp\Big\{x\in\R^m\Big|\,\sum_{i=1}^m x_i=1,\;x\ge 0 \Big\}
\end{eqnarray*}
describes the set of mixed strategies for the $x$-player
(Player~1) with $m$ pure strategies; similarly for the $y$-player
$y\in\D_n$ (Player~2). This means that if Player~1 uses $x\in
\D_m$ and Player~2 uses $y\in \D_n$, then Player~1 gets payoff
$-x\transp A y$ while Player~2 gets payoff $x\transp A y$. Thus
the equilibrium problem \eqref{eq:equil} can be reformulated as
the following problem of {\em nonsmooth convex optimization}:
\begin{eqnarray}\label{opt}
\mbox{minimize }\;F(x,y)\;\mbox{ subject to }\;(x,y)\in\DD,
\end{eqnarray}
where the minimizing cost function $F(x,y)$ is defined by the
maximum
\begin{eqnarray}\label{eq:DefF}
F(x,y):=\max\big\{x\transp A v-u\transp A
y\big|\;(u,v)\in\DD\big\}.
\end{eqnarray}
It is easy to observe that
\begin{eqnarray}\label{opt1}
\min\{F(x,y)|\;(x,y)\in\DD\}=0.
\end{eqnarray}
Taking \eqref{opt1} into account, we say \cite{GilPenaSandh,s}
that a feasible pair $(\ox,\oy)\in\DD$ is a {\em Nash equilibrium}
to \eqref{eq:equil} if $F(\ox,\oy)=0$, which corresponds to an
optimal solution of the constrained optimization problem
\eqref{opt}. Consider the optimal solution set
\begin{eqnarray}\label{Nash}
S:=\big\{(\ox,\oy)\in\DD\big|\;F(\ox,\oy)=0\big\}=F^{-1}(0)\cap
(\DD)
\end{eqnarray}
and, following \cite{GilPenaSandh}, define the {\em condition
measure} $\k(A)$ of the matrix game \eqref{eq:equil} depending on
the underlying matrix $A$ via the objective \eqref{eq:DefF} and
the optimal solution set \eqref{Nash} as
\begin{equation}\label{eq:k}
\k(A):=\inf\left\{\k\ge 0\,\bigl|\,\dist\big((x,y);S\big)\le\k
F(x,y)\;\mbox{ for all }\;(x,y)\in\DD\right\},
\end{equation}
where $\dist(\cdot;S)$ stands for the standard Euclidean distance
function.\vspace*{0.05in}

In what follows we derive three major results concerning the
characterization of the condition measure $\k(A)$ in \eqref{eq:k}.
The  first theorem shows that the condition measure $\k(A)$
precisely relates to the exact bound of metric regularity for an
associated set-valued mapping built upon the cost function
\eqref{eq:DefF}. The second result expresses this exact regularity
bound via the subdifferential of the convex function
\eqref{eq:DefF} and the normal cone to the simplex product $\DD$
and then computes the latter constructions in terms of the initial
data of \eqref{eq:equil}. Finally, we arrive at an exact formula
for evaluating $\k(A)$, which is a key step towards performing
further complexity analysis of the algorithm
\cite{GilPenaSandh}.\vspace*{0.05in}

To formulate the first theorem, define a set-valued mapping
$\Phi\colon\R^{m+n}\tto\R$ by
\begin{eqnarray}\label{Phi}
\Phi(x,y):=\left\{\begin{array}{ll}
\big[F(x,y),\infty)\;&\mbox{ if }\;(x,y)\in\DD,\\\\
\emp\;&\mbox{ otherwise}
\end{array}
\right.
\end{eqnarray}
via the cost function $F$ constructed in \eqref{eq:DefF}. Let
$\reg\Phi\left((x,y),F(x,y)\right)$ be the {\em exact bound of
metric regularity} (or the {\em exact regularity bound/modulus})
of the mapping $\Phi$ around the point
$\left((x,y),F(x,y)\right)\in\gph\Phi$; see
\cite{Morduk2006Book1,RockWets1998}. For the reader's convenience
we recall these concepts in Section~2 below.

\begin{theorem}\label{th1}{\bf (condition measure via the exact
regularity bound).} Assume that $(\DD)\setminus S\ne\emp$ with $S$
defined in \eqref{Nash}. Then we have the precise relationship
\begin{eqnarray}\label{k-reg}
\k(A)=\disp\sup_{(x,y)\in(\DD)\setminus
S}\reg\Phi\big((x,y),F(x,y)\big)
\end{eqnarray}
between the condition measure \eqref{eq:k} and the exact
regularity bound of \eqref{Phi}.
\end{theorem}

Our second major result gives a complete
characterization of metric regularity for the set-valued mapping
$\Phi$ defined in \eqref{Phi}.

\begin{theorem}\label{th2} {\bf (computing the exact bound of metric regularity).}
For any point $(x,y)\in(\DD)\setminus S$, the exact regularity
bound of the mapping $\Phi$ from \eqref{Phi} around the point
$((x,y),F(x,y))$ admits the representation
\begin{eqnarray}\label{regPhi}
\reg\Phi\big((x,y),F(x,y)\big)=\disp\frac{1}{\dist\big(0;\partial
F(x,y)+N_\DD(x,y)\big)}
\end{eqnarray}
via the subdifferential of the convex function \eqref{eq:DefF} and
the normal cone to the simplex product $\DD$ at $(x,y)$.
\end{theorem}

Unifying the results of Theorem~\ref{th1} and Theorem~\ref{th2}
and then explicitly computing the subdifferential and normal
cone on the right-hand side of \eqref{regPhi}, we get the precise
formula for computing the condition measure of the smoothing
algorithm for matrix games as stated in Theorem~\ref{th3} below.
To formulate this third major result, we introduce some convenient
notation. Let $a_i$ as $i=1,\ldots,n$ and $-b\transp_k$ as
$k=1,\ldots,m$ stand for the columns and the rows of the matrix
$A$, respectively. By  $e_j$, $j=1,\ldots,m+n$, we denote the unit
vectors in $\R^{m+n}$, i.e.,
\begin{eqnarray*}
(e_j)_l=0\;\mbox{ for all }\;l\ne j\;\mbox{ and }\;(e_j)_j=
1\;\mbox{ as }\;j=1,\ldots,m+n.
\end{eqnarray*}
For a positive integer $p$, let ${\bf 1}_p:=\matr{1&\dots& 1}\in
\R^p$.  Finally, given a feasible point
$(x,y)\in\Delta_m\times \Delta_n$, define the index sets
$I(x), K(y),$ and $J(x,y)$ by
\begin{eqnarray}\label{ind}
\left\{\begin{array}{ll} I(x):=\Big\{\bar \imath\in
\{1,\ldots,n\}\Big|\;a_{\bar \imath}\transp
x=\disp\max_{i\in\{1,\ldots,n\}}a_i\transp x\Big\},\\ K(y):=
\Big\{\bar k\in\{1,\ldots,m\}\Big|\;b_{\bar k}\transp y
=\disp\max_{k\in\{1,\ldots,m\}}b_k\transp y\Big\},\\
J(x,y):=\big\{j\in\{1,\ldots,m\}\big|\;x_j=0\big\}\bigcup\big\{j=m+p\big|\;y_p
=0\big\}.
\end{array}\right.
\end{eqnarray}

\begin{theorem}\label{th3} {\bf (computing the condition
measure).} Let $(\DD)\setminus S\ne\emp$. Then, in the notation
above, the condition measure $k(A)$ defined in \eqref{eq:k} is
computed by
\begin{eqnarray*}
\begin{array}{ll}
\k(A)=\disp\sup_{(x,y)\in(\DD)\setminus
S}\Big[&\dist\Big(0;\co\big\{(a_i,b_k)\big|\;i\in I(x),\;k\in
K(y)\big\}\\
&+\span\{{\bf 1}_m\}\times\span\{{\bf 1}_n\}
-\cone\big[\co\big\{e_j\big|\;j\in
J(x,y)\big\}\big]\Big)\Big]^{-1},
\end{array}
\end{eqnarray*}
where the symbols {\rm span}, {\rm cone},
and {\rm co} stand respectively for the linear, conic, and convex
hulls of the sets in question.
\end{theorem}

The proofs of Theorem~\ref{th1} and Theorem~\ref{th2} given below
are based on applying advanced techniques of variational analysis
and generalized differentiation. This approach leads us therefore
to deriving the precise formula for the condition measure in
Theorem~\ref{th3}. For additional insight, we also present a
direct, independent proof of the latter theorem relying on more
conventional while somewhat more laborious techniques of convex
optimization employing particularly Lagrangian duality.

\begin{remark}\label{numer} {\bf (numerical implementation and
further research).} {\rm Numerical implementation of the formula
for the condition measure in Theorem~\ref{th3} is not a purpose of
this paper and in fact is not an easy job. It has been well
recognized in complexity theory that evaluating condition measures
may be in general as difficult as to solve the original problem.
This is true, e.g., in the cases of such fundamental complexity
measures as the condition number of a matrix \cite{Higham96} used
in estimating complexity of numerical linear algebra algorithms,
Renegar's condition number \cite{PeRe,Renegar95a,Renegar95b} that
characterizes difficulty of solving conic feasibility problems,
the ``measure of condition" for finding zeros of complex
polynomials introduced by Shub and Smale \cite{ShubSmaleBezoutI},
etc.

The main purpose of this paper is not obtaining an easily
computable expression for the condition measure $\k(A)$, but
rather gaining a better understanding on how exactly the problem
data influence the condition measure. Observe that the formula for
$\k(A)$ obtained in Theorem~\ref{th3} is much easier to evaluate
and analyze than the original construction (\ref{eq:k}). This is
valuable for the average-case and smoothed analysis of the
algorithm, singling out classes of well-conditioned problems,
preconditioning issues, and making further improvements to the
algorithm. We will pursue these goals in our subsequent research.}
\end{remark}

The rest of the paper is organized as follows. In
Section~\ref{sec:MRGD} we recall some basic definitions and facts
of variational analysis and generalized differentiation crucial
for deriving the main results of the paper.
Section~\ref{sec:PfMainTheorem} is devoted to variational proofs
of the main results formulated above. Finally, in
Section~\ref{sec:PfExactK} we present an alternative direct proof
of Theorem~\ref{th3} by employing tools of convex optimization.

Throughout the paper we use standard notation and terminology of
variational analysis; see, e.g., the basic texts
\cite{Morduk2006Book1,RockWets1998}.

\section{Preliminaries from variational analysis and generalized
differentiation}\label{sec:MRGD}
\setcounter{equation}{0}

Here we confine ourselves to finite-dimensional Euclidean spaces
sufficient for the subsequent applications. The reader is referred
to \cite{Morduk2006Book1,RockWets1998} for more details and
related material.

Given a set-valued mapping $G\colon\R^n\tto\R^m$, consider its
inverse $G^{-1}\colon\R^m\tto\R^n$ with $x\in G^{-1}(z)
\Longleftrightarrow z\in G(x)$ as well as its graphs, domain, and
range defined respectively by
\begin{eqnarray*}
\gph G:=\big\{(x,z)\big|\;z\in G(x)\big\},\qquad\dom
G:=\big\{x\big|\;G(x)\ne\emp\big\},\qquad\rge G:=\dom G^{-1}.
\end{eqnarray*}
The notion of metric regularity is of primary interest in our
development.

\begin{definition}\label{mr} {\bf (metric regularity).} A
set-valued mapping $G\colon\R^n\tto\R^m$ is {\sc metrically
regular} around $(\ox,\oz)\in\gph G$ with modulus $\mu\ge 0$ if
there exist neighborhoods $U$ of $\ox$ and $V$ of $\oz$ such that
\begin{eqnarray}\label{mr1}
\dist\big(x;G^{-1}(z)\big)\le\mu\dist\big(z;G(x)\big)\;\mbox{
whenever }\;x\in U\;\mbox{ and }\;z\in V.
\end{eqnarray}
The infimum of $\mu\ge 0$ over all $(\mu,U,V)$ for which
\eqref{mr1} holds is called the {\sc exact regularity bound} of
$G$ around $(\ox,\oz)$ and is denoted by $\reg G(\ox,\oz)$.
\end{definition}

It is well known in variational analysis that the fundamental
property of metric regularity is closely related to Lipschitzian
behavior of inverse mappings. Recall that a mapping
$G\colon\R^n\tto\R^m$ is {\em Lipschitz-like} (or has the Aubin
property) around $(\ox,\oz)\in\gph G$ with modulus $\ell\ge 0$ if
there are neighborhoods $U$ of $\ox$ and $V$ of $\oz$ such that
\begin{eqnarray}\label{lip}
G(x)\cap V\subset G(u)+\ell\|x-u\|\B\;\mbox{ for all }\;x,u\in U,
\end{eqnarray}
where $\B$ stands for the Euclidean closed unit ball of the space
in question.
The infimum of $\ell\ge 0$ over all the combinations $(\ell,U,V)$
for which \eqref{lip} holds is called the {\em exact Lipschitzian
bound} of $G$ around $(\ox,\oz)$ and is denoted by $\lip
G(\ox,\oz)$.\vspace*{0.05in}

The following result can be found, e.g., in
\cite[Theorem~1.49]{Morduk2006Book1}.

\begin{proposition}\label{thm:RegLip} {\bf (relationships between
metric regularity and Lipschitz-like properties).} Let
$G\colon\R^n\tto\R^m$, and let $(\ox,\oz)\in\gph G$. Then the
mapping $G$ is metrically regular around $(\ox,\oz)$ if and only
if its inverse $G^{-1}\colon\R^m\tto\R^n$ is Lipschitz-like around
$(\oz,\ox)$. Furthermore, we have the equality
\begin{eqnarray*}
\reg G(\ox,\oz)=\lip G^{-1}(\oz,\ox).
\end{eqnarray*}
\end{proposition}

One of the key advantages of modern variational analysis is the
possibility to completely characterize Lipschitzian and metric
regularity properties of set-valued mappings in terms of
appropriate generalized differential constructions enjoying full
calculus. Let us recall such constructions used in this paper.

Given a nonempty subset $\O\subset\R^n$ and a point $\ox\in\O$,
define the Fr\'echet/regular normal cone to $\O$ at $\ox$ by
\begin{eqnarray*}
\Hat N_\Omega(\ox):=\left\{v\in\R^n\,\Bigr|\,
\limsup_{x^{\underrightarrow{\Omega}}\ox}\frac{\langle v,
x-\ox\rangle}{\|x-\ox\|}\le 0\right\},
\end{eqnarray*}
where the symbol $x\st{\O}{\to}\ox$ means that $x\to\ox$ with
$x\in\O$. Then the Mordukhovich (basic/limiting) {\em normal cone}
to $\O$ at $\ox\in\O$ is defined by
\begin{eqnarray}\label{nc}
N_\O(\ox):=\disp\Limsup_{x\st{\O}{\to}\ox}\Hat N_\O(x),
\end{eqnarray}
where `$\Limsup$' stands for the Painlev\'e-Kuratowski outer/upper
limit of a set-valued mapping $M\colon\R^n\tto\R^m$ given by
\begin{eqnarray*}
\begin{array}{ll}
\disp\Limsup_{x\to\ox}M(x):=\Big\{v\in\R^m\Big|&\exists\,x_k\to\ox,\;v_k\to
v\;\mbox{ as }\;k\to\infty\;\mbox{ such that }\\
&v_k\in M(x_k)\;\mbox{ for all }\;k=\N:=\{1,2,\ldots\}\Big\}.
\end{array}
\end{eqnarray*}
If the set $\O$ is locally closed around $\ox$, the normal cone
\eqref{nc} admits the equivalent description (which was in fact
the original definition in \cite{M76})
\begin{eqnarray*}
N_\O(\ox)=\disp\Limsup_{x\to\ox}\big[\cone\big(x-\Pi_\O(x)\big)\big]
\end{eqnarray*}
in terms of the projection operator
$\Pi_\O(x):=\{y\in\O\,|\;\|y-x\|=\dist(x;\O)\}$.

Note that the normal cone \eqref{nc} may be nonconvex even for
simple sets $\O\subset\R^n$, e.g., for the graph of $|x|$ and the
epigraph of $-|x|$ at $(0,0)$. Due to its nonconvexity, the normal
cone \eqref{nc} cannot be polar to any tangent cone. Nevertheless,
this nonconvex normal cone and the corresponding
subdifferential/coderivative constructions for
extended-real-valued (i.e., with values in $(-\infty,\infty]$)
functions and set-valued mappings satisfy comprehensive calculus
rules, which are derived by using variational arguments,
particularly the extremal principle of variational analysis; see
\cite{Morduk2006Book1,RockWets1998} and the references therein.

A set $\O$ is called {\em normally regular} at $\ox\in\O$ if
$N_\O(\ox)=\Hat N_\O(\ox)$. The class of normally regular sets
covers ``nice" sets having a local convex-like structure. A major
example is provided by convex sets; see, e.g.,
\cite[Proposition~1.5]{Morduk2006Book1}.

\begin{proposition}\label{prop:ConvexRegular} {\bf (normal regularity
of convex sets).} Let $\O\subset\R^n$ be convex. Then it is
normally regular at every point $\ox\in\O$, and its normal cone
\eqref{nc} reduces to the normal cone in the sense of convex
analysis:
\begin{eqnarray*}
N_\O(\ox)=\big\{v\in\R^n\big|\;\la v,x-\ox\ra\le 0\;\mbox{ for all
}\;x\in\O\big\}.
\end{eqnarray*}
\end{proposition}

Given next a set-valued mapping $G\colon\R^n\tto\R^m$ and a point
$(\ox,\oz)\in\gph F$, define a generalized derivative of $G$ at
$(\ox,\oz)$ induced by the normal cone \eqref{nc} to $\gph G$ at $(\ox,\oz)$.
 Namely, the {\em coderivative} of $G$ at $(\ox,\oz)$
is a set-valued mapping $D^*G(\ox,\oz)\colon\R^m\tto\R^n$ with the
values
\begin{eqnarray}\label{cod}
D^*G(\ox,\oz)(v):=\big\{u\in\R^n\big|\;(u,-v)\in N_{\gph
G}(\ox,\oz)\big\}.
\end{eqnarray}

Observe that $0\in D^*G(\ox,\oz)(0)$ and $D^*G(\ox,\oz)(\lm v)=\lm
D^* G(\ox,\oz)(v)$ for every $\lambda>0$, i.e., the coderivative \eqref{cod} is a
positively homogeneous mapping. If $G\colon\R^n\to\R^m$ is
single-valued and smooth around $\ox$ with the derivative $\nabla
G(\ox)$, we have (see, e.g.,  by \cite[Theorem 1.38]{Morduk2006Book1})
\begin{eqnarray*}
D^*G(\ox)(v)=\big\{\nabla G(\ox)\transp v\big\}\;\mbox{ for all
}\;v\in\R^m.
\end{eqnarray*}
The latter signifies that the coderivative \eqref{cod} is an
appropriate extension of the adjoint/transpose derivative operator
to the case of nonsmooth and set-valued mappings. Note also that,
by the nonconvexity of the normal cone \eqref{nc}, the
coderivative \eqref{cod} is not dual to any tangentially generated
graphical derivative, except of the case when $G$ is {\em
graphical regular} at $(\ox,\oz)$ meaning that
\begin{eqnarray*}
N_{\gph G}(\ox,\oz)=\Hat N_{\gph G}(\ox,\oz).
\end{eqnarray*}

As mentioned above, the coderivative \eqref{cod} satisfies
comprehensive calculus rules for general set-valued mappings. In
this paper we only need the following one, which is a consequence
of \cite[Proposition~3.12]{Morduk2006Book1}. To formulate it,
recall that the indicator mapping $\dd_{\O}\colon\R^n\to\R$ of a
set $\O\subset\R^n$ is defined by
\begin{eqnarray*}
\dd_\O(x):=\left\{\begin{array}{ll} 0&\mbox{if }\;x\in\O,\\
\emp &\mbox{otherwise},
\end{array}\right.
\end{eqnarray*}
(a bit different from the indicator functions) and that we easily
have the relationship
\begin{eqnarray*}
D^*(\delta_{\Omega})(\ox)(v)=N_\Omega (\ox)\;\mbox{ for any
}\;\ox\in\O\;\mbox{ and }\;v\in\R.
\end{eqnarray*}
\begin{proposition}\label{prop:CodSum} {\bf (coderivative sum rule).}
Let $\O\subset\R^n$ be locally closed around $\ox\in\O$, and let
$G\colon\R^n\rightrightarrows\R$ be closed-graph and
Lipschitz-like around $(\ox,\oz)\in\gph G$. Then for all $v\in\R$
we have the inclusion
\begin{eqnarray*}
D^*(G+\dd_\O)(\ox,\oz)(v)\subset D^*G(\ox,\oz)(v)+N_\O(\ox),
\end{eqnarray*}
which holds as equality if $\O$ is normally regular at $\ox$ and
$F$ is graphically regular at $(\ox,\oz)$.
\end{proposition}

In what follows we employ the norm of the coderivative as a
positively homogeneous mapping. The norm of a positively
homogeneous mapping $M\colon\R^n\tto\R^m$ is defined by
\begin{eqnarray*}
\|M\|:=\sup\big\{\|u\|\,\big|\,u\in M(v)\;\mbox{ with }\;\|v\|\le
1\big\}
\end{eqnarray*}
and admits (by passing to the inverse) the useful distance
function representation below established in
\cite[Proposition~2.5]{DonLewRock}.

\begin{proposition}\label{prop:norm} {\bf (norm of positively homogeneous
mappings).} Let $M\colon\R^n\tto\R^m$ be positively homogeneous.
Then the norm of its inverse is computed by
\begin{eqnarray*}
\|M^{-1}\|=\sup_{\|v\|=1}\frac{1}{\dist\big(0;M(v)\big)}.
\end{eqnarray*}
\end{proposition}

The final and most important result presented in this section
provides a complete coderivative characterization of the
Lipschitz-like property (known as the Mordukhovich criterion
\cite{RockWets1998}) with computing the exact bound of
Lipschitzian moduli; see \cite[Theorem~5.7]{M93},
\cite[Theorem~4.10]{Morduk2006Book1}, and
\cite[Theorem~9.40]{RockWets1998} for different proofs.

\begin{theorem}\label{thm:CodLip}{\bf (coderivative characterization of the
Lipschitz-like property for set-valued mappings).} Let
$G\colon\R^m\tto\R^n$ be closed-graph around $(\bar x,\bar z)\in
\gph G$. Then $G$ is Lipschitz-like around this point if and only
if $D^*G(\ox,\oz)(0)=\{0\}$. In this case
\begin{eqnarray*}
\lip G(\bar x,\bar z)=\|D^*G(\bar x,\bar z)\|.
\end{eqnarray*}
\end{theorem}

\section{Proofs of main results}\label{sec:PfMainTheorem}
\setcounter{equation}{0}

We give here complete proofs of Theorem~\ref{th1} and
Theorem~\ref{th2}, and thus derive the condition measure formula
of Theorem~\ref{th3} by variational arguments.

Let us start with the proof of Theorem~\ref{th1}. To proceed, we
first establish a more convenient representation of the condition
measure \eqref{eq:k} for our further analysis, which in turn is
preceded by a technical claim.

Observe that the function $F(x,y)$ defined by (\ref{eq:DefF}) can
be written as follows:
$$
F(x,y) =\max_{\begin{array}{c} i=1,\dots,n\\  j = 1, \dots,
m\end{array}} (a_i\transp x+b_j\transp y).
$$
In addition we represent the simplex product $\DD$ by:
$$
\DD = \big\{w =(x,y)\, |\, w\geq 0,\; Ew=f\big\} \quad
\mbox{with}\quad E: = \left[
\begin{array}{cc}
{\bf 1}_m\transp & 0 \\
0 & {\bf 1}_n\transp
\end{array}\right]\;\mbox{ and }\;f: =  \left[
\begin{array}{c}
1 \\
1
\end{array}\right].
$$
To simplify notation, rewrite the function $F$ as
$$
F(w) = \max_{\ell\in L}c_\ell \transp w,
$$
where $L: =\{1,\dots,n\} \times \{1,\dots,m\}$ and $c_\ell\transp:
= \matr{a_i\transp & b_j\transp}$ for each $\ell = (i,j) \in L$.
Denote further ${\Omega}:= \Delta_n \times \Delta_m = \{w\,|\,
w\geq 0, Ew = f\}$ and rewrite \eqref{eq:k} as
\begin{equation}\label{eq:SimplK}
k_F:= \inf\big\{k\ge 0\big|\;\dist (w; S) \leq k F(w)\;\mbox{ for
all }\;w\in{\Omega}\big\}
\end{equation}
with $S$ given by \eqref{Nash}. Observe that $\min_{\Omega} F(x) =
0$ by \eqref{opt1}. It is also convenient for us to define the
moving sets
\begin{eqnarray}\label{S}
S(z) := \big\{w\in {\Omega}\big|\; F(w) = z\big\}=F^{-1}(z)\cap
{\Omega}\;\mbox{ with thus }\;S=S(0)
\end{eqnarray}
and to represent the mapping $\Phi$ in \eqref{Phi} and its inverse
by
\begin{equation}\label{eq:PhiInv}
\Phi(w)= \big[F(w),\infty\big) +
\delta_{\Omega}(w)\quad\mbox{and}\quad \Phi^{-1}(z) = \big\{w\in
{\Omega}\big|\; F(w)\leq z\big\}.
\end{equation}
Let us finally denote $\J:=\{1,\dots,n,n+1,\dots,m+n\}$ and define
the corresponding counterparts of  the index sets  $I(\cdot)$ and
$J(\cdot)$ from \eqref{ind} given by
\begin{eqnarray*}
\I(w):=\big\{\ell \in L\big|\; c_\ell \transp w = F(w)\big\},\quad
\J(w):=\big\{j\in \J\big|\;w_j = 0\big\}.
\end{eqnarray*}

It is not difficult to verify the following technical claim, where
$\B_\gg$ stands for the closed ball of radius $\gg>0$ centered at
the origin.

\begin{claim}\label{lem:IndexNeigh} {\bf (relationships between index sets).}
For every $\bar w \in {\Omega}$ there exists $\gg>0$ such that
$\I(w)\subset \I(\bar w)$ and $\J(w)\subset \J(\bar w)$ whenever
$w\in \bar w+\B_\gg$.
\end{claim}
{\bf Proof.} Fix an arbitrary element $\bar w\in\Omega$ and let
$$
0<\gg < \min\left\{\min_{\begin{array}{c}\ell\in L,\\ c_\ell \ne
0\end{array}}\frac{1}{2\|c_\ell\|} \min_{i\in L\setminus \I(\bar
w)}(F(\bar w) - c_i\transp \bar w), \min_{j\in \J\setminus
\J(\bar w)}{\bar w_j}\right\},
$$
where $\min\emp=\infty$ by convention. It is easy to observe that
such a number $\gg$ always exists. When $\I(\bar w) = L$, the
inclusion $\I(w)\subset \I(\bar w)$ is obvious. Assume thus that
$L\setminus\I(\bar w)\ne \emptyset$. For every $\ell_0\in
L\setminus \I(\bar w)$ and every $w\in \bar w +\B_\gg$ we have
$$
c_{\ell_0}\transp w  =   c_{\ell_0}\transp \bar w +
c_{\ell_0}\transp (w- \bar w) \leq c_{\ell_0}\transp \bar w +
\|c_{\ell_0}\|\gg < c_{\ell_0}\transp\bar w + \frac{1}{2}
\min_{\ell\in L\setminus \I(\bar w)}\big(F(\bar w) - c_\ell\transp
\bar w\big),
$$
which implies the relationships
\begin{equation}\label{eq:005}
\begin{array}{ll}
\disp\max_{\ell\in L\setminus\I(\bar w)}( c_{\ell}\transp w)&<
\disp\max_{\ell\in L\setminus\I(\bar w)}(c_{\ell}\transp\bar w) +
\frac{1}{2} \min_{\ell\in L\setminus \I(\bar w)}\big(F(\bar w) -
c_\ell\transp
\bar w\big)\\\\
& = \disp\frac{1}{2}\left( F(\bar w)+ \max_{i\in L\setminus
\I(\bar w)}c\transp_\ell \bar w\right).
\end{array}
\end{equation}
Similarly, for every $\ell_0 \in \I(\bar w)$ we have
$$
c_{\ell_0}\transp w  =  c_{\ell_0}\transp \bar w +
c_{\ell_0}\transp (w- \bar w) \ge  c_{\ell_0}\transp \bar w -
\|c_{\ell_0}\|\gg > c_{\ell_0}\transp\bar w  - \frac{1}{2}
\min_{\ell\in L\setminus\I(\bar w)}\big(F(\bar w) - c_\ell\transp
\bar w\big),
$$
which in turn implies that
\begin{equation}\label{eq:006}
\begin{array}{ll}
\disp\max_{\ell\in \I(\bar w)}( c_{\ell}\transp w
)&>\disp\max_{\ell\in \I(\bar w)}(c_\ell\transp\bar w) -
\disp\frac{1}{2} \min_{\ell\in
L\setminus \I(\bar w)}\big(F(\bar w) - \disp c_\ell\transp \bar w \big)\\\\
& =\disp\max_{\ell\in \I(\bar w)}(c_\ell\transp\bar w)
-\disp\frac{1}{2} F(\bar w)+ \frac{1}{2}\max_{\ell\in L\setminus
\I(\bar w)}\disp (c_\ell\transp \bar w)
\\\\
& = \disp\frac{1}{2}\left( F(\bar w)+ \max_{\ell\in L\setminus
\I(\bar w)}\disp (c_\ell\transp \bar w)\right),
\end{array}
\end{equation}
where the last step follows by the construction of $\I(\bar w)$.
Combining (\ref{eq:005}) and (\ref{eq:006}) gives us the strict
inequality
$$
\max_{\ell\in \I(\bar w)}( c_\ell\transp w)- \max_{\ell\in
L\setminus\I(\bar w)}( c_\ell\transp w)> 0,
$$
and hence justifies the first claimed inclusion $\I(w)\subset
\I(\bar w)$.

It remains to show that $\J(w)\subset\J(\bar w)$ for all
$w\in\ow+\B_{\gg}$. The latter inclusion is obvious when $\J(\bar
w)=\J$. Assume now that $\J \setminus \J(\bar w)\ne \emptyset$ and
then get for every $w\in \bar w+ \B_\gg$ and $j\in \J\setminus
\J(\bar w)$ the relationships
$$
w_j = \bar w_j +w_j - \bar w_j > \bar w_j -\gg \geq \bar w_j -
\min_{j\in \J\setminus \J(\bar w)} {\bar w_j} \geq 0.
$$
Thus $w_j > 0$ whenever $j\in \J\setminus \J(\bar w)$, and thus we
arrive at $\J(w)\subset \J(\bar w)$. $\h$\vspace*{0.05in}

The next result provides a useful representation of the condition
measure \eqref{eq:SimplK} convenient for our subsequent analysis.

\begin{lemma}\label{lem:SupAttained} {\bf (representation of
condition measure).} Assume that ${\Omega}\setminus S \ne
\emptyset$. Then there exists $\bar w\in {\Omega}\setminus S$ such
that
\begin{equation}\label{eq:supAttained}
k_F = \sup_{w\in {\Omega}\setminus S} \frac{\dist(w; S)}{F(w)} =
\frac{\dist (\bar w;S)}{F(\bar w) }.
\end{equation}
\end{lemma}

The proof of Lemma~\ref{lem:SupAttained} is based on dividing the
set $\Omega\setminus S$ into a finite family of subsets and
showing that the supremum is attained on each one of such sets.
Before presenting this proof, we introduce some notation and prove
a technical proposition.

For every $s\in S$ let $N_S(s)$ be the normal cone to $S$ at $s$. Observe that this cone is
polyhedral and can be represented as follows:
\begin{equation}\label{eq:016}
N_S(s) =\big\{g=w-s \,|\, \dist(w;S)=\|w-s\|\big\}.
\end{equation}
Also for every $s\in S$ define
$$
K_s: =\Omega \cap\big(s+N_S(s)\big).
$$
Note that for every $s\in S$ the set $K_s$ is a convex polytope
and that
$$
\Omega = \bigcup_{s\in S}K_s \quad \mbox{and}\quad \Omega\setminus
S = \bigcup_{s\in S}\big(K_s\setminus \{s\}\big).
$$
In addition, it follows from \eqref{eq:016} and the definition of
$K_s$ that for every $s\in S$ and any $w\in K_s$ we have
$$
\dist(w;S)=\|w-s\|.
$$

\begin{proposition}\label{prop:Ks} {\bf (supremum attainability).}
For every $s\in S$ there exists $\bar w\in K_s\setminus \{s\}$,
which realizes the supremum
\begin{equation}\label{eq:012}
\sup_{w\in K_s\setminus\{s\}}\frac{\dist(w;S)}{F(w)}
=
\frac{\dist(\bar w;S)}{F(\bar w)}.
\end{equation}
\end{proposition}

\noindent
{\bf Proof. } Let $s\in S$, and let a sequence $\{w_k\}\subset
K_s\setminus \{s\}$ be such that
$$
\sup_{w\in K_s\setminus\{s\}}\frac{\dist(w;S)}{F(w)} =\sup_{w\in K_s
\setminus\{s\}}\frac{\|w-s\|}{F(w)} =
\lim_{k\to\infty }\frac{\|w_k-s\|}{F(w_k)}.
$$
Without loss of generality it can be assumed that $w_k\to w_0\in
K_s$ along the whole sequence $\{w_k\}$, since $K_s$ is closed and
bounded. If $w_0\neq s$, we get $F(w_0)\neq 0$ and thus
$$
\lim_{k\to \infty} \frac{\|w_k-s\|}{F(w_k)} = \frac{\|w_0-s\|}{F(w_0)}.
$$
Therefore, in this case we have found $\bar w = w_0$  that
satisfies \eqref{eq:012} .

In the case when $w_0 =s$, consider the sequence $\{g_k\}$ defined
by
$$
g_k: = \frac{w_k-s}{\|w_k-s\|}\;\mbox{ for all }\;
k\in\N.
$$
Since $\{g_k\}$ is bounded, suppose without loss of generality
that $g_k\to g$ with $\|g\|=1$. The polyhedrality of $K_s$ implies
the existence of $\lambda_1>0$ satisfying $s+\lambda g\in K_s$ for
all $\lambda \in[0,\lambda_1]$. Further, the finiteness of the
index set $L$ allows us to get without loss of generality that
$$
{\cal I}(w_k) = {\cal I}_0 \subset L\;\mbox{ for all}\;k.
$$
Then there exists $l_0\in {\cal I}_0$ such that
\begin{equation}\label{eq:010}
c_{l_0}\transp w_k \geq c_l\transp w_k\;\mbox{ whenever }\;l\in
L\;\mbox{ and }\;k\in \N.
\end{equation}
Since $w_k \to s$, we have for all $l\in L$ that
$$
\lim_{k\to \infty} c_l\transp w_k =c_l\transp s,
$$
which by \eqref{eq:010} yields that $l_0\in {\cal I}(s)$.

Let us show in addition that there is $\lambda_2>0$ for which
$l_0\in {\cal I}(s+\lambda g)$ whenever $\lambda\in
[0,\lambda_2]$. Indeed, assume the contrary and find a sequence
$\{\lambda^{(m)}\}$ with $\lambda^{(m)}\downarrow 0$ such that
$l_0\notin {\cal I}(s+\lambda^{(m)} g)$ for all $m\in \N$. Then
without loss of generality (since $L$ is finite) there exists
$\bar l \in L$ satisfying
$$
c_{\bar l}\transp \big(s+\lambda^{(m)} g\big) >
c_{l_0}\transp\big(s+\lambda^{(m)} g\big)\;\mbox{ for all }\;m\in
\N.
$$
By the same arguments used above to show that $l_0\in {\cal
I}(s)$, we get $\bar l \in {\cal I}(s)$ and thus
\begin{equation}\label{eq:011}
c_{l_0}\transp s = c_{\bar l}\transp s = 0.
\end{equation}
Together with \eqref{eq:010}, the latter yields that
$$
c_{l_0}\transp\frac{(w_k-s)}{\|w_k-s\|}\geq c_{\bar l}\transp\frac{(w_k-s)}{\|w_k-s\|}.
$$
By passing to the limit as $k\to \infty$, we obtain
$c_{l_0}\transp g \geq c_{\bar l}\transp g$ and conclude that
$$
c_{l_0}\transp\big(s+\lambda^{(m)}g\big) \ge c_{\bar l}\transp
\big(s+\lambda^{(m)}g\big),
$$
which contradicts our assumption and thus justifies the claim.

To proceed further, denote $\lambda:=\min\{\lambda_1,\lambda_2\}$
and $\bar w: = s+\lambda g $. Then we have
$$
\frac{\dist(\bar w;S)}{F(\bar w)} = \frac{\|\bar w -s \|}{F(\bar
w)} =
\frac{\lambda \|g\|}{\lambda c_{l_0}\transp g} =
\frac{1}{c_{l_0}\transp g}.
$$
But it follows at the same time that
$$
\lim_{k\to \infty}\frac{\|w_k -s \|}{F(w_k)} =
\lim_{k\to\infty}\frac{\|w_k - s\|}{ c_{l_0}\transp
w_k}=\lim_{k\to \infty}\frac{1}{
c_{l_0}\transp\frac{w_k-s}{\|w_k-s\|}} = \frac{1}{c_{l_0}\transp
g},
$$
which shows that $\bar w$ satisfies \eqref{eq:012} and thus
completes the proof of the proposition. $\h$\vspace*{0.1in}

Now let us justify our basic
Lemma~\ref{lem:SupAttained}.\vspace*{0.05in}

\noindent {\bf Proof of Lemma~\ref{lem:SupAttained}.} For every
$s\in S$ consider the tangent cone to $\Omega$ at $s$ defined by
$$
T_{\Omega}(s): =\big\{g\in \R^n\big|\, \exists\, \alpha_0>0\,:\,
s+\alpha g \in \Omega \;\mbox{ for all }\;\alpha\in
(0,\alpha_0]\big\}.
$$
We split the set $S$ into a family of disjoint subsets ${\cal S}$,
which correspond to the following equivalence classes:
$$
s_1\sim s_2\quad \mbox{if and only if}\quad  N_S(s_1) = N_S(s_2),
\quad T_\Omega(s_1) =  T_\Omega(s_2),\quad\mbox{and}\quad
\partial F(s_1) = \partial F(s_2).
$$
Observe the relationship
\begin{equation}\label{sup}
\sup_{w\in \Omega \setminus S}\frac{d(w,S)}{F(w)}=\sup_{C\in {\cal
S}}\sup_{s\in C}\sup_{w\in K_s \setminus
\{s\}}\frac{d(w,S)}{F(w)},
\end{equation}
where the outer supremum (over $C\in {\cal S}$) on the right-hand
side is attained, since the number of different sets/classes in
$\cal S$ is finite due to the polyhedral structure of the problem.
The innermost supremum in \eqref{sup} is attained by
Proposition~\ref{prop:Ks}. Hence to prove the claim, it remains to
show that the supremum over $s\in C$ in \eqref{sup} is also
attained. We do it by proving that for every $C\in {\cal S}$ it
holds
\begin{equation}\label{eq:015}
\sup_{w\in K_{s_1} \setminus
\{s_1\}}\frac{d(w,S)}{F(w)}=\sup_{w\in K_{s_2} \setminus
\{s_2\}}\frac{d(w,S)}{F(w)}\;\mbox{ whenever }\;s_1,s_2\in C.
\end{equation}

To proceed, pick an arbitrary set $C\in{\cal S}$ and arbitrary
elements $s_1,s_2\in C$. Fix a point $w_0 \in K_{s_1}\setminus
\{s_1\}$ and get by the definition of $K_{s_1}$ the corresponding
point $g:=w_0-s_1\in N_{S}(s_1)$. Since $N_{S}(s_1) = N_{S}(s_2)$
and since $N_{S}(s_2)$ is a cone, we have
$$
\lambda g \in N_{S}(s_2)\;\mbox{ for all }\;\lambda \in
[0,\infty).
$$
It follows again from the definition of $K_s$ and from the
polyhedrality of $\Omega$ that $g\in T_\Omega(s_1)$ and hence
$g\in T_\Omega(s_2)$ by $T_\Omega(s_1)=T_\Omega(s_2)$. This
ensures the existence of $\lambda_1>0$ such that $s_2+\lambda g
\in \Omega$ for all $\lambda\in [0,\lambda_1]$. Thus
$$
s_2+\lambda g \in \Omega \cap\big(s_2+N_S(s_2)\big)
\setminus\{s_2\} = K_{s_2} \setminus\{s_2\}\;\mbox{ for all
}\;\lambda \in (0,\lambda_1].
$$
Since $F$ is polyhedral, it follows that for each point $s$ in $S$
the function $F(\cdot)-F(s)$ is positively homogeneous in a
neighborhood of $s$, and therefore it can be represented in this
neighborhood via the subdifferential of $F$ as
\begin{equation*}
F(w)-F(s)=\disp\max_{v\in\partial F(s)}v\transp(w-s).
\end{equation*}
Since $\partial F(s_1) = \partial F(s_2)$, find $\lambda_2>0$ such
that $F(s_1+\lambda g)= F(s_2+\lambda g)$ for all $\lambda \in
[0,\lambda_2]$. Letting further $\lambda: = \min\{\lambda_1,
\lambda_2\}$ and denoting $w_1:=s_1+\lambda g\in K_{s_1}$ and
$w_2: = s_2+\lambda g\in K_{s_2}$, we have $F(w_1)=F(w_2)$ and
$$
F(w_1) = F(s_1+\lambda g) = F\big((1-\lambda)s_1+\lambda
w_0\big)\le (1-\lambda)F(s_1) +\lambda F(w_0) = \lambda F(w_0)
$$
by the convexity of $F$. This implies the relationships
$$
\frac{\dist(w_2;S)}{F(w_2)} =\frac{\dist(w_1;S)}{F(w_1)}=
\frac{\|\lambda g\|}{F(w_1)} \geq \frac{\lambda \|g\|}
{\lambda F(w_0)}=\frac{\dist(w_0;S)}{F(w_0)}.
$$
Since the latter holds for any $w_0 \in K_{s_1}\setminus \{s_1\}$,
we get
\begin{equation}\label{sup1}
\sup_{w\in K_{s_2}\setminus \{s_2\}}\frac{\dist(w;S)}{F(w)} \ge
\sup_{w\in K_{s_1}\setminus \{s_1\}}\frac{\dist(w;S)} {F(w)}.
\end{equation}
The inverse inequality to \eqref{sup1} is obtained by
interchanging the roles of $s_1$ and $s_2$. Since our initial
choice of $C\in {\cal S}$ and $s_1,s_2\in C$ was arbitrary, we
arrive at the equality in \eqref{eq:015} for all $C\in {\cal S}$
and thus complete the proof of the lemma. $\h$\vspace*{0.1in}

Now we are ready to prove our main results, namely
Theorem~\ref{th1}, Theorem~\ref{th2} and Theorem~\ref{th3}.

\subsection{Proof of Theorem~\ref{th1}}
We split the proof of the theorem into three major
steps.\vspace*{0.1in}

\noindent {\bf Step~1: metric regularity via condition measure.}
{\em For every $\bar w\notin S$ and every $\bar z >0$ we have the
distance estimate}
\begin{equation}\label{eq:claimDistPhi}
\dist\big(\bar w; \Phi^{-1}(\bar z)\big)\leq k_F
\dist\big(\Phi(\bar w);\bar z\big).
\end{equation}

\noindent {\bf Proof.} When $\bar w\notin {\Omega}$, the right-hand side of
(\ref{eq:claimDistPhi}) becomes infinity (by the construction of
$\Phi$ in \eqref{eq:PhiInv} and the standard convention on
$\inf\emp=\infty$) while the left-hand side is finite, i.e., there
is nothing to prove. Considering the case of $\bar w\in {\Omega}$,
observe that the left-hand side of (\ref{eq:claimDistPhi}) becomes
zero when $F(\bar w)\le\bar z$, and thus the inequality holds
automatically. It remains to examine the case when $0<\bar
z<F(\bar w)$ with $\bar w\in\O$.

To proceed, let $w^*:=\Pi_{\Phi^{-1}(\bar z)}(\bar w)$, and observe
that $F(w^*)=\bar z$, since
otherwise the continuity of $F$ would allow us to find a closer point to $w$
in $[w^*,\bar w]\cap \Phi^{-1}(\bar z)$. Thus
\begin{equation}\label{eq:claim101}
\dist\big(\bar w; \Phi^{-1}(\bar z)\big) = \dist \big(\bar w;
S(\bar z)\big) =\|\bar w - w^*\|,
\end{equation}
where $S(\cdot)$ is defined in \eqref{S}. Let $w_0:=\Pi_S(\bar
w).$ Employing again the continuity of  $F$, we find $\lambda \in
(0,1)$ such that
\begin{equation}\label{eq:claim102}
F(w_0+\lambda (\bar w - w_0)) = \bar z.
\end{equation}
In addition the convexity of $F$ yields that
\begin{equation}\label{eq:claim103}
\bar z = F\big(w_0 +\lambda (\bar w - w_0)\big)\leq F(w_0)
+\lambda \big(F(\bar w)- F(w_0)\big)= \lambda F(\bar w).
\end{equation}
Combining the above, we have the relationships
\begin{eqnarray}\label{eq:claim104}
\begin{array}{ll}
\dist\big(\bar w; \Phi^{-1}(\bar z)\big)&= \dist \big(\bar
w; S(\bar z)\big)\qquad\big(\mbox{ by \eqref{eq:claim101}}\big)\\
&\le  \|w_0 +\lambda (\bar w - w_0)- \bar w\|
\qquad\big(\mbox{ by \eqref{eq:claim102}}\big)\\
&=(1-\lambda)\|w_0 - \bar w\| = (1-\lambda)\dist (\bar w; S);
\end{array}
\end{eqnarray}
\begin{eqnarray}\label{eq:claim105}
\begin{array}{ll}
\dist\big(\Phi(\bar w); \bar z\big) & =F(\bar w) - \bar z
\qquad\big(\mbox{ as $\bar z<F(\bar w)$}\big)\\
& \ge  F(\bar w) - \lambda F(\bar w)
\qquad \big(\mbox{ by \eqref{eq:claim103}}\big)\\
& =(1-\lambda) F(\bar w),
\end{array}
\end{eqnarray}
which finally give
\begin{eqnarray*}
\begin{array}{ll}
\dist\big(\bar w; \Phi^{-1}(\bar z)\big)) & \leq (1-\lambda)\dist
(\bar w; S)\qquad \big(\mbox{ by \eqref{eq:claim104}}\big)\\
& \leq    (1-\lambda)k_F F(\bar w)
\qquad\big(\mbox{ as $\dist(\bar w;S)\le k_F F(\bar w)$}\big)\\
& \leq   k_F\dist \big(\Phi(\bar w); \bar z\big) \qquad
\big(\mbox{ by \eqref{eq:claim105}}\big)
\end{array}
\end{eqnarray*}
and thus allow us to arrive at \eqref{eq:claimDistPhi}.
$\h$\vspace*{0.1in}

\noindent{\bf Step~2: distance properties.} {\em Let $\bar w\in
{\Omega}\setminus S$ be such that
\begin{equation}\label{slope}
k_F = \frac{\dist(\bar w; S)}{F(\bar w)},
\end{equation}
let $w_0:=\Pi_S(\bar w)$, and for $\lambda
\in (0,1)$ let
$$
w_\lambda : = \bar w +\lambda (w_0-\bar w).
$$
Then for any $\lm\in(0,1)$ we have the properties}:
\begin{itemize}
\item[{\bf (i)}] $F(w_\lambda) = (1-\lambda) F(\bar w)$.
\item[{\bf (ii)}] $\dist(\bar w; \Phi^{-1}(F(w_\lambda)))=\lambda
\dist(\bar w; S)$.
\end{itemize}

\noindent {\bf Proof.} To justify (i), observe that by the convexity of $F$
we have
\begin{equation}\label{eq:claim201}
F(w_\lambda)\leq F(\bar w) +\lambda \big(F(w_0)-F(\bar w)\big)
= (1-\lambda) F(\bar w)
\end{equation}
in the notation above. On the other hand, the definition of $k_F$
and the choice of $\bar w$ yield
\begin{equation}\label{eq:claim202}
F(w_\lambda)\geq \frac{\dist(w_\lambda;S)}{k_F} =
\frac{(1-\lambda)\|\bar w - w_0\|}{k_F}=
\frac{(1-\lambda)\dist(\bar w; S)}{k_F}= (1-\lambda) F(\bar w).
\end{equation}
Thus assertion (i) follows from \eqref{eq:claim201} and
\eqref{eq:claim202}.\vspace*{0.05in}

\noindent To justify (ii), it suffices to show that
$$
\dist\big(\bar w; \Phi^{-1}(F(w_\lambda))\big) = \|\bar
w-w_\lambda\|.
$$
Proceeding by contradiction, assume that $\dist\big(\bar w;
\Phi^{-1}(F(w_\lambda))\big) < \lambda\dist(\bar w;S)=\lambda k_F
F(\bar w)$ and let $w^*:=\Pi_{S(F(w_\lambda))}(\bar w)$. By the
continuity of $F$ we have
$$
\dist \big(\bar w; \Phi^{-1}(F(w_\lambda))\big)=\dist\big(\bar w;
S(F(w_\lambda))\big)=\| \bar w - w^*\|,
$$
which yields therefore that
\begin{equation}\label{claim.eqn1}
\|\bar w- w^* \| = \dist\big(\bar w; \Phi^{-1}(F(w_\lambda))\big)
< \lambda  k_F  F(\bar w).
\end{equation}
Taking further a point $\tilde w$ closest to $w^*$ in $S$, we get
by \eqref{eq:SimplK} and by part (i) above that
\begin{equation}\label{claim.eqn2}
\|\tilde w - w^*\| \le k_F F(w^*) = k_F F(w_\lambda) = k_F
(1-\lambda)F(\bar w).
\end{equation}
Since $\tilde w\in S$, the latter implies that
\begin{eqnarray*}
\begin{array}{ll}
\dist(\bar w;S) &\le\|\bar w -  \tilde w\| \\
&\leq \|\bar w - w^*\| + \|w^* -  \tilde w\| \qquad \text{(by the
triangle inequality)}\\
&<\lambda k_F F(\bar w) +(1-\lambda)k_F F(\bar w) \qquad
\big(\text{by \eqref{claim.eqn1}
and \eqref{claim.eqn2}}\big) \\
&=k_F F(\bar w),
\end{array}
\end{eqnarray*}
which contradicts \eqref{slope} and thus completes the proof of
Step~2. $\h$\vspace*{0.1in}

\noindent{\bf Step~3: condition measure via metric regularity.}
{\em We have the equality}
\begin{eqnarray}\label{k1}
k_F = \sup_{w\in \Omega \setminus S} \reg \Phi \big(w, F(w)\big).
\end{eqnarray}

\noindent {\bf Proof.} Let us first show that
\begin{eqnarray}\label{k2}
k_F\ge\sup_{w\in {\Omega}\setminus S}\reg\Phi\big(w,F(w)\big).
\end{eqnarray}
Assuming the contrary, find $(w', z')\in \gph F$, $w'\in
{\Omega}\setminus S$ satisfying
$$
\reg \Phi(w', z') > k_F.
$$
Observe that there exists a neighborhood of $(w', z')$ such that
for all points $(w,z)$ in that neighborhood $w\notin S$ and $z>0$.
By the definition of metric regularity we can find $\bar w, \bar
z$ in such a neighborhood of $(w', z')$ for which
$$
\dist\big(\bar w, \Phi^{-1}(\bar z)\big)>k_F\dist\big(\Phi(\bar
w); \bar z\big).
$$
The latter contradicts Step~1 and thus ensures \eqref{k2}.
\vspace*{0.05in}

To prove the opposite inequality in \eqref{k1}, by
Lemma~\ref{lem:SupAttained} find $\bar w\in\O\setminus S$ such
that
$$
\dist(\bar w; S) = k_F F(\bar w).
$$
Let $w_0:=\Pi_S(\bar w)$ and define
$$
w_\lambda: = \bar w+ \lambda(w_0 - \bar w),\quad 0<\lm<1.
$$
It follows from Step~2 that for every $\lambda\in (0,1)$ and the
above choice of $\bar w$ we have
$$
\frac{\dist\big(\bar w;
\Phi^{-1}(F(w_\lambda))\big)}{\dist(F(w_\lambda); \Phi\big(\bar
w)\big)} = \frac{\lambda\dist(\bar w;S)}{\lambda F(\bar w)} = k_F.
$$
The latter implies, since $w_\lambda \to \bar w$ and
$F(w_\lambda)\to F(\bar w)$ as $\lambda \downarrow 0$, that
$$
\reg \Phi\big(\bar w, F(\bar w)\big) \geq \limsup_{\lambda
\downarrow 0} \frac{\dist\big(\bar w; \Phi^{-1}(F
(w_\lambda))\big)}{\dist\big(F(w_\lambda); \Phi(\bar w)\big)} =
k_F,
$$
which therefore yields
$$
k_F \leq \sup_{w\in {\Omega}\setminus S} \reg \Phi\big(w,F(w)\big)
$$
and completes the the proof of the theorem. $\h$

\subsection{Proof of Theorem~\ref{th2}}

First of all, observe by Step~1 in the proof of Theorem~\ref{th1}
that the multifunction $\Phi$ is metrically regular around $(w,
z)\in \gph \Phi$ for every $w\in {\Omega}\setminus S$. Employing
the corresponding results of Section~\ref{sec:MRGD}, for $(w,z)
\in \gph \Phi$ with $w\in{\Omega}\setminus S$ we get
\begin{eqnarray*}
\begin{array}{ll}
\reg \Phi(w,z) &= \lip \Phi^{-1}(z,w)\qquad \big(\mbox{ by Theorem~\ref{thm:RegLip}}\big)\\
&= \|D^* \Phi^{-1}(z,w)\|\qquad \big(\mbox{ by Theorem~\ref{thm:CodLip}}\big) \\
&= \big\|\big(D^* \Phi(w,z)\big)^{-1}\big\|\qquad \big(\mbox{ by the definition of
$D^*\Phi(w, z$)}\big)
\\\\
&= \disp\sup_{|v|=1}\frac{1}{ \dist\big(0; D^*
\Phi(w,z)(v)\big)}\qquad \big(\mbox{ by
Proposition~\ref{prop:norm} with $n=1$}\big).
\end{array}
\end{eqnarray*}
This gives therefore the regularity exact bound formula
\begin{equation}\label{eq:reg2}
\reg \Phi(w,z) = \frac{1}{\min\big\{\dist\big(0; D^*
\Phi(w,z)(-1)\big),\dist\big(0; D^* \Phi(w,z)(1)\big)\big\}}.
\end{equation}

Defined next a set-valued mapping $\Tilde
F\colon\R^{m+n}\rightrightarrows \R$ by
\begin{eqnarray}\label{tilde F}
\Tilde F(w):= \big[F(w),\infty\big)\;\mbox{ with }\;\gph\Tilde F =
\epi F
\end{eqnarray}
and observe that it is Lipschitz-like at every point of its graph,
which is the epigraph of a Lipschitz continuous function.
Furthermore, the graph of $\Tilde F$ is convex, and hence $\Tilde
F$ is graphically regular at any point of its graph by
Proposition~\ref{prop:ConvexRegular}, which also ensures the
normal regularity of the convex set $\O$. Applying
Proposition~\ref{prop:CodSum} to the sum $\Phi=\Tilde F+\dd_\O$,
we get the equality
\begin{equation}\label{eq:Step201}
D^*\Phi(w, z)(\lambda) = D^*\Tilde F(w,z)(\lambda) +N_{\Omega}(w)\;\mbox{
for all }\;\lm\in\R.
\end{equation}
It follows from the structure of $\Tilde F$ in \eqref{tilde F},
the coderivative definition \eqref{cod}, and the well-known
subdifferential representation
\begin{eqnarray*}
\partial\ph(\ox)=\big\{v\in\R^n\big|\;(v,-1)\in
N\big((\ox,\ph(\ox));\epi\ph\big)\big\},\quad\ox\in\dom\ph,
\end{eqnarray*}
for any convex function $\ph\colon\R^n\to\oR$ that
\begin{equation}\label{eq:Step202}
D^*\Tilde F(w, z)(1) = \left\{
\begin{array}{ll}
\dc F(w), & z= F(w),\\
\emptyset, & z> F(w),
\end{array}
\right. \qquad D^*\Tilde F(w, z)(-1) = \emptyset.
\end{equation}
Combining \eqref{eq:reg2}, \eqref{eq:Step201}, and
\eqref{eq:Step202} gives us the formula
\begin{eqnarray}\label{reg1}
\reg \Phi(w,z)=\frac{1}{\dist\big(0;\partial F
(w)+N_\Omega(w)\big)},
\end{eqnarray}
which is (\ref{regPhi}). $\h$

\subsection{Proof of Theorem~\ref{th3}.}
By Theorem~\ref{th1} and Theorem~\ref{th2}, it suffices to show
the following representations for $\partial F(x,y)$ and
$N_{\Delta_m \times \Delta_n}(x,y)$:
\begin{eqnarray}\label{subF}
\partial F(x,y)=\co\big\{(a_i,b_k)\in\R^m\times\R^n\big|\;i\in
I(x),\;k\in K(y)\big\},
\end{eqnarray}
\begin{eqnarray}\label{normalD}
N_{\DD}(x,y)= \span\{{\bf 1}_m\}\times\span\{{\bf
1}_n\}-\cone\big[\co\big\{e_j\big|\;j\in J(x,y)\big\}\big].
\end{eqnarray}
Indeed,
the classical subdifferential formula for max-functions (see,
e.g., \cite[Exercise~8.31]{RockWets1998}) gives us
\begin{eqnarray*}
\dc\big(\max_{\ell\in L}c_\ell\transp w\big)  = \co\big\{c_{\bar
\ell}\big|\; \bar \ell \in L, \; c_{\bar \ell}\transp w =
\max_{\ell\in L} c_\ell\transp w \big\}.
\end{eqnarray*}
This implies by the max-structure of the function $F$ in
\eqref{eq:DefF} that
\begin{eqnarray*}
\dc F(x, y)  = \co\big\{(a_i, b_k)\big|\;i\in I(x),\;k\in
K(y)\big\},
\end{eqnarray*}
which is \eqref{subF}.

To prove~\eqref{normalD}, we recall first the calculus formula
$$
N_{A\cap B}(w) = N_A(w)+N_B(w)
$$
held at every $w\in A\cap B$ for arbitrary convex polyhedra in
finite dimensions; see, e.g., \cite[Corollary~23.8.1]{R-CA}). Thus
we have in our case that
$$
N_{\Omega}(w)= N_{\{Eu= f\}}(w) + N_{\{u\geq 0\}}(w).
$$
Moreover, it is easy to see that
$$
N_{\{Eu = f\}}(w) = ({\ker}E)^\perp = \span \{{\bf 1}_m\}\times
\span \{{\bf 1}_n\},
$$
\begin{eqnarray*}
N_{\{u\ge 0\}}(w) = -\cone\big[\co\big\{e_j\, |\, j\in
\J(w)\big\}\big].
\end{eqnarray*}
Thus for any $w = (x,y) \in \Omega = \Delta_n \times \Delta_m$ we
have
$$
N_{\Delta_m \times \Delta_n}(x,y) = \span \{{\bf 1}_m\}\times\span
\{{\bf 1}_n\} - \cone\big[\co\big\{e_j\, |\, j\in
J(x,y)\big\}\big],
$$
which is \eqref{normalD}. $\h$

\section{Condition measure formula via alternative proof}\label{sec:PfExactK}
\setcounter{equation}{0}

In this section we give another proof of Theorem~\ref{th3} based
on convex optimization. This proof is split into three lemmas and
the preceding technical claim.\vspace*{0.05in}

Given a point $\bar w\in {\Omega}\setminus S$ and keeping the
notation above, consider the following two problems of parametric
optimization (with the parameter $z\in\R$) defined by
$$
\begin{array}{ccc}
\begin{array}{l}
\quad V_z(\bar w): = \min\limits_w  \|w-\bar w\|\\
\begin{array}{crcl}
\;\mbox{ s.t.} & c_\ell\transp w & \leq & z \quad \forall \ell\in \I(\bar w)\\
& Ew &=& f\\
& w_j & \geq & 0 \quad \forall j \in \J(\bar w)
\end{array}\quad (P_z)
\end{array}
& \mbox{and} &
\begin{array}{l}
\quad \Tilde V_z(\bar w): = \min\limits_w  \|w-\bar w\|\\
\begin{array}{crcl}
\;\mbox{ s.t.} & \max\limits_{\ell\in L}\{c_\ell\transp w\} & = &  z\\
& Ew &=& f\\
& w & \geq & 0
\end{array}\quad  (\Tilde P_z)
\end{array}
\end{array}
$$
and name $(P_z)$ and $(\Tilde P_z)$ the {\em first} and {\em
second parametric problem}, respectively. Observe that for every
$\bar w \in {\Omega}\setminus S$ and $z\in \R_+$ the optimal value
$\Tilde V_z(\bar w)$ in problem ($\Tilde P_z$) is equal to $\dist
(\bar w; S(z))$. Although the proof of the following claim is
straightforward, we provide it for completeness and the reader's
convenience.

\begin{claim}\label{claim:DeltaEpsilon} {\bf (stability of optimal solutions
to first parametric problem).} For any $\bar w \in {\Omega}
\setminus S$ and any $\gg > 0$ there is $\e>0$ $($depending on
$\bar w$ and $\gg)$ such that whenever $z\in [F(\bar w) - \e,
F(\bar w)]$ a unique solution $w_z$ to problem $(P_z)$ exists and
satisfies the continuity property $\|w_z - \bar w\|\leq \gg$ with
respect to the parameter $z$.
\end{claim}

\noindent {\bf Proof.} Fix $\gg>0$ and put $w^S:=\Pi_S(\bar w)$. Let
$$
\tilde w : = \bar w +\tau (w^S - \bar w)\;\mbox{ with }\;  \tau: =
\min\left\{\frac{\gg}{\|w^S- \bar w\|}, 1\right\}\in (0,1].
$$
Setting $\e:= F(\bar w) - F(\tilde w)$, observe by the convexity
of $F$ that
$$
\e =  F(\bar w) - F(\tilde w)  \geq  \tau F(\bar w) > 0.
$$
We have furthermore that
$$
E \tilde w =  E \bar w +  \tau  (Ew^S - E\bar w) = f,
$$
$$
\tilde w_j = (1-\tau)\bar w_j + \tau w^S_j\geq 0\;\mbox{ for all
}\; j\in \J(\bar w),\;\mbox{ and}
$$
$$
\|\tilde w - \bar w\| =   \tau  \|w^S - \bar w\|\leq \gg,
$$
which imply the inclusion
$$
\tilde w \in \Delta := \big\{ w\big|\; \|w-\bar w\|\leq \gg,\;\,
Ew = f,\;\, w_j\geq 0 \;\mbox{ for all }\;j\in\J(\bar w)\big\}.
$$
Since the set $\Delta$ is obviously convex with $\bar w\in
\Delta$, we get
$$
w_t: = \bar w + t (\tilde w - \bar w) \in \Delta\;\mbox{ whenever
}\; t \in [0,1].
$$
It follows from  $F(\tilde w)= F(\bar w) - \e$ and the continuity
of $F$ that for every $z\in [F(\bar w) - \e, F(\bar w)]$ there is
$t_z \in [0,1]$ such that $w_{t_z}$ satisfies the equation
$$
F(w_{t_z})= z.
$$
For any $z$ from the above we easily get that $w_{t_z}$ is
feasible to problem ($P_z$), that the set of feasible solutions to
this problem is surely closed and bounded, and that the cost
function is continuous with respect to $w$. Thus $(P_z)$ admits an
optimal solution, which is unique as a unique projection of $\bar
w$ on the convex feasible set. Finally,
$$
V_z(\bar w) = \|w_z- \bar w\|\leq \|w_{t_z} - \bar w\| \leq t_z
\|\tilde w - \bar w\|\leq \gg,
$$
and hence the optimal solution $w_z$ belongs to the ball $\bar
w+\B_{\gg}$. $\h$\vspace*{0.05in}

The next result, whose proof is based on
Claims~\ref{lem:IndexNeigh} and~\ref{claim:DeltaEpsilon},
indicates the parameter region on which the optimal values in the
first and second parametric problems agree.

\begin{lemma}\label{lem:OptCoincide} {\bf (optimal values agree
for both parametric problems).} Let $\bar w\in {\Omega}\setminus
S$. Then there exists $\e_{\bar w}\in(0,F(\bar w))$ such that for
every parameter $z\in [F(\bar w)-\e_{\bar w},F(\bar w)]$ the
optimal values of problems $(P_z)$ and $(\Tilde P_z)$ coincide.
\end{lemma}

\noindent
{\bf Proof.} Fix $\bar w\in {\Omega}\setminus S$ and observe that
the set of feasible solutions for ($\Tilde P_z$) obviously belongs
to the set of feasible solutions for $(P_x)$. Thus we have $\Tilde
V_z(\bar w)\ge V_z(\bar w)$ for all $z\in\R$. It remains to show
that there exists $\e_{\bar w}>0$ such that $\Tilde V_z(\bar w)\le
V_z(\bar w)$ whenever $z\in [F(\bar w)-\e_{\bar w},F(\bar w)]$.

Employing Claim~\ref{lem:IndexNeigh}, find $\gg>0$ for which
$\I(w) \subset \I(\bar w)$ and $\J(w) \subset \J(\bar w)$ when
$w\in \bar w +\B_\gg$. Further, it follows from
Claim~\ref{claim:DeltaEpsilon} that for such $\gg$ there is $\e>0$
with the property: whenever $z\in [F(\bar w) - \e , F(\bar w)]$
there exists a unique solution $w_z$ to problem (${ P}_z$)
satisfying $w_z\in \bar w+ \B_\gg$. Our choice of $\gg$ ensures
the feasibility of $w_z$ in problem ($\Tilde P_z$), and therefore
we have the relationships
$$
\Tilde V_z(\bar w) \leq \|w_z-\bar w\| =  V_z(\bar w),
$$
which thus complete the proof of the lemma. $\h$

\begin{lemma}\label{lem:GlobIsLoc} {\bf (distances to solution sets).}
For every $ w\in {\Omega}\setminus S$ denote $z_{ w}: = F( w)-
\e_{ w}$ with $\e_{w}$ taken from
Lemma~{\rm\ref{lem:OptCoincide}}. Then we have
\begin{eqnarray}\label{sup-k}
\sup_{w\in {\Omega}\setminus S} \frac{\dist\big( w;
S(z_{w})\big)}{F(w)- z_{w}} = \sup_{w\in {\Omega}\setminus S}
\frac{\dist(w; S)}{F(w)}.
\end{eqnarray}
\end{lemma}

\noindent
{\bf Proof.} Fix $w\in {\Omega}\setminus S$ and $z\in (0,F(w))$,
and then let
$$
k( w, z):= \frac{\dist(w; S(z))}{F(w)-z}.
$$
We first justify the inequality
\begin{equation}\label{eq:kzxleqk}
\sup_{ w\in {\Omega}\setminus S} k( w, z_{w})\leq \sup_{w \in
{\Omega}\setminus S} k(w,0),
\end{equation}
which gives the one in \eqref{sup-k}. Pick any $\bar w \in
{\Omega} \setminus S$, let $w^S: = \arg\min_{w\in S}\|w-\bar w\|$
and
$$
w_t: = (1-t)w^S+t\bar w,\quad t\in [0,1].
$$
Since $z_{\bar w}\in (0, F(\bar w))$, there is $\tau\in (0,1)$
with $F(w_\tau) = z_{\bar w}$. By the convexity of $F$ we have
\begin{equation}\label{eq:zxleqtF}
z_{\bar w} = F(w_\tau) \leq (1-\tau)F(w^S)+\tau F(\bar w) = \tau F(\bar w).
\end{equation}
Further, it follows from $w_\tau \in S(z_{\bar w})$ that
\begin{equation}\label{eq:dist1}
\dist \big(\bar w; S(z_{\bar w})\big) \leq \|\bar w -
w_\tau\|=(1-\tau)\|w^S- \bar w\| = (1-\tau) \dist (\bar w; S).
\end{equation}
Combining (\ref{eq:zxleqtF}) and (\ref{eq:dist1}) gives us
$$
k(\bar w, z_{\bar w}) = \frac{\dist\big(\bar w; S(z_{\bar
w})\big)}{\dist\big(F(\bar w); z_{\bar w}\big)}\leq \frac{(1-\tau)
\dist(\bar w;S)}{(1-\tau)F(\bar w)} = k(\bar w, 0),
$$
which yields (\ref{eq:kzxleqk}) and thus the corresponding
inequality in \eqref{sup-k} .

It remains to show that
\begin{equation}\label{eq:kzxgeqk}
\sup_{w\in {\Omega}\setminus S} k( w, z_w)\ge\sup_{w \in
{\Omega}\setminus S} k( w, 0),
\end{equation}
which ensures the equality in \eqref{sup-k}. By
Lemma~\ref{lem:SupAttained} we have that $k_F = \sup_{w\in
{\Omega}\setminus S} k(w, 0)$ and that the maximum is attained at
some $\bar w\in {\Omega}\setminus S$. Given $z \in (0,F(\bar w))$,
let
$$
w_z := {\rm arg}\min_{w\in S(z_w)}\|w-\bar w\|
$$
and observe the estimate
$$
\dist(\bar w; S)\leq \|\bar w -  w_z\|+\dist ( w_z; S),
$$
implying in turn that
$$
\sup_{w\in {\Omega}\setminus S} k( w, z_w)\geq k(\bar w, z_{\bar
w})=
\frac{\|\bar w - w_{z_{\bar w}}\|} {F(\bar w) - z_{\bar w}}
\geq
 \frac{\dist(\bar w; S)-\dist ( w_{z_{\bar w}}; S)}{F(\bar w) - z_{\bar w}}.
$$
On the other hand, we have the equality $\dist(\bar w; S) = F(\bar
w)k(\bar w,0) $ by the definition of $k(w,z)$ and also the
relationships
\begin{eqnarray*}
\dist(w_{z_{\bar w}};S) = F(w_{z_{\bar w}}) k(w_{z_{\bar w}},0)\le
z_{\bar w} k(\bar w,0)
\end{eqnarray*}
due to $F(w_{z_{\bar w}}) = z_{\bar w}$ and $k(\bar w,0) =
\sup_{w\in \Omega\setminus S} k(w,0)$. Thus
$$
\sup_{w\in {\Omega}\setminus S} k( w, z_w) \geq \frac{F(\bar
w)k(\bar w, 0)-z_{\bar w} k(\bar w, 0)}{F(\bar w) - z_{\bar w}}=
k(\bar w, 0) = \sup_{w \in
 {\Omega}\setminus S} k( w,0),
$$
which justifies \eqref{eq:kzxgeqk} and completes the proof of the
lemma. $\h$\vspace*{0.05in}

The last lemma establishes, by employing Lagrangian duality, a
precise formula for computing the optimal value of the cost
function in the parametric problem $(P_z)$---and hence in $(\Tilde
P_z)$---via the initial data.

\begin{lemma}\label{lem:slnSubprob} {\bf (computing optimal values of
parametric problems).} Let $\bar w \in {\Omega} \setminus S$ and
$z\in (0,F(\bar w))$. Then the optimal value $V_z(\bar w)$ of
problem $(P_z)$ is computed by
\begin{eqnarray*}
V_z(\bar w) = \disp\frac{F(\bar w) - z}{\dist\big(0;
\co\big\{c_i, i\in \I(\bar w)\big\}+(\ker
E)^\perp-\cone\big[\co\big\{e_j,\, j\in \J(\bar w)\big\}\big]\big)}.
\end{eqnarray*}
\end{lemma}

\noindent
{\bf Proof.} Observe that problem ($P_z$) can be reformulated as
$$
V_z(\bar w) = \inf_w \sup_{
\begin{array}{c}
\|u\|\leq 1,\\
\lambda_i \geq 0, i\in \I(\bar w),\\
v\in \R^m,\\
\mu_j \geq 0, j \in \J(\bar w),\\
\mu_j = 0, j \in \J \setminus \J(\bar w)
\end{array}
                         }
\left[ u\transp (w - \bar w) +\sum_{i\in \I(\bar w)}\lambda_i
(c\transp_i w - z) +v\transp(Ew-f)-\sum_{j\in \J(\bar w)}\mu_j
w_j\right].
$$
Observe that the convex optimization problem $(P_z)$ satisfies the
Slater condition: For $\delta \in (0,1)$ sufficiently
small and $\tilde w \in S$, the point $(1-\delta) \tilde w +
\delta (\frac{1}{n} {\bf 1}_n, \frac{1}{m} {\bf 1}_m)$ is a
strictly feasible point. Therefore, we can interchange the
supremum and the infimum above by Lagrangian duality. This gives
$$
V_z(\bar w) = \sup_{
\begin{array}{c}
\|u\|\leq 1,\\
\lambda_i \geq 0, i\in \I(\bar w),\\
v\in \R^m,\\
\mu_j \geq 0, j \in \J(\bar w),\\
\mu_j = 0, j \in \J \setminus \J(\bar w)
\end{array}
} \inf_w \left[ u\transp (w - \bar w) +\sum_{i\in \I(\bar
w)}\lambda_i (c\transp_i w - z) +v\transp(Ew-f)-\sum_{j\in \J(\bar
w)}\mu_j w_j\right].
$$
Regrouping the terms inside the square brackets, we obtain
$$
V_z(\bar w) =   \sup_{
\begin{array}{c}
\|u\|\leq 1,\\
\lambda_i \geq 0, i\in \I(\bar w),\\
v\in \R^m,\\
\mu_j \geq 0, j \in \J(\bar w),\\
\mu_j = 0, j \in \J \setminus \J(\bar w)
\end{array}
} \inf_w \left[ \left(u+\sum_{i\in \I(\bar
w)}\lambda_ic_i+E\transp v-\mu\right)\transp w - u\transp\bar w
-\sum_{i\in \I(\bar w)}\lambda_i z - v\transp f \right].
$$
Observe further that, whenever the term $(u+\sum_{i\in \I(\bar
w)}\lambda_ic_i+E\transp v-\mu) $ is not zero, the inner infimum
in $w$ necessarily becomes $-\infty$. This allows us to put
$$
u=-\sum_{i\in \I(\bar w)}\lambda_ic_i- E\transp v+\mu
$$
and consequently rewrite the expression for $V_z(\bar w)$ as
follows:
$$
V_z(\bar w) =    \sup_{\begin{array}{c}
\|\sum_{i\in \I(\bar w)}\lambda_i c_i +E\transp v - \mu\|\leq 1,\\
\lambda_i \geq 0, i\in \I(\bar w),\\
\mu_j \geq 0, j \in \J(\bar w),\\
\mu_j = 0, j \in \J \setminus \J(\bar w)
\end{array}}
\left[\left(\sum_{i\in \I(\bar w)}\lambda_i c_i+E\transp
v-\mu\right)\transp \bar w -\sum_{i\in \I(\bar w)}\lambda_i z -
v\transp f \right].
$$
Regrouping again gives us the formula
$$
V_z(\bar w) =           \sup_{
\begin{array}{c}
\|\sum_{i\in \I(\bar w)}\lambda_i c_i +E\transp v - \mu\|\leq 1,\\
\lambda_i \geq 0, i\in \I(\bar w),\\
\mu_j \geq 0, j \in \J(\bar w),\\
\mu_j = 0, j \in \J \setminus \J(\bar w)
\end{array}
}\left[\sum_{i\in \I(\bar w)}\lambda_i\left(c_i\transp \bar w
-z\right)+v\transp(E\bar w- f) -\mu\transp \bar w \right].
$$
Noting that $E\bar w= f$, $\mu\transp \bar w = 0$ and $c_i\transp
\bar w = F(\bar w)$ for $i\in \I(\bar w)$, we have
$$
V_z(\bar w) =           \left(F(\bar w) - z\right) \sup_{
\begin{array}{c}
\|\sum_{i\in \I(\bar w)}\lambda_i c_i +E\transp v - \mu\|\leq 1,\\
\lambda_i \geq 0, i\in \I(\bar w),\\
\mu_j \geq 0, j \in \J(\bar w),\\
\mu_j = 0, j \in \J \setminus \J(\bar w)
\end{array}
} \sum_{i\in \I(\bar w)}\lambda_i.
$$
Since $V_z(\bar w)\ne 0$, the latter yields
$$
V_z(\bar w) =          \left(F(\bar w) - z\right) \sup_{
\begin{array}{c}
\left\|\frac{\sum_{i\in \I(\bar w)}\lambda_i c_i}{\sum_{i\in
\I(\bar w)}\lambda_i} +E\transp \frac{v}{\sum_{i\in \I(\bar
w)}\lambda_i} - \frac{\mu}{\sum_{i\in \I(\bar w)}
\lambda_i}\right\|\leq \frac{1}{  \sum_{i\in \I(\bar w)}\lambda_i},\\
\lambda_i \geq 0, i\in \I(\bar w),\\
\mu_j \geq 0, j \in \J(\bar w),\\
\mu_j = 0, j \in \J \setminus \J(\bar w)
\end{array}
} \frac{1}{1/\sum_{i\in\I(\bar w)}\lambda_i},
$$
which can be written as
$$
V_z(\bar w) =           \left(F(\bar w) - z\right)\sup_{
\begin{array}{c}
\lambda_i \geq 0, i\in \I(\bar w),\\
\mu_j \geq 0, j \in \J(\bar w),\\
\mu_j = 0, j \in \J\setminus \J(\bar w)
\end{array}
} \frac{1}{\left\|\frac{\sum_{i\in \I(\bar w)}\lambda_i
c_i}{\sum_{i\in \I(\bar w)}\lambda_i} +E\transp
\frac{v}{\sum_{i\in \I(\bar w)}\lambda_i} - \frac{\mu}{\sum_{i\in
\I(\bar w)}\lambda_i}\right\|}.
$$
Changing further the variables by
$$\tilde \lambda_i: =
\disp\frac{\lambda_i}{\sum_{i\in \I (\bar w)}\lm_i},\quad \tilde
\mu:= \disp\frac{\mu}{\sum_{i\in \I (\bar w)}\lm_i},\quad \tilde v
:=\disp\frac{v}{\sum_{i\in \I (\bar w)}\lm_i},
$$
we arrive at the expression
$$
V_z(\bar w)    =          \left(F(\bar w) - z\right) \sup_{
\begin{array}{c}
\tilde \lambda_i \geq 0, i\in \I(\bar w),\\
\sum_{i\in \I(\bar w)}\tilde \lambda_i = 1,\\
\tilde \mu_j \geq 0, j \in \J(\bar w),\\
\tilde \mu_j = 0, j \in \J\setminus \J(\bar w)
\end{array}
}\frac{1}{\left\|\sum_{i\in \I(\bar w)}\tilde \lambda_i c_i
+E\transp \tilde v - \tilde \mu\right\|},
$$
which can be equivalently written as
$$
V_z(\bar w) =      \frac{F(\bar w) - z}{\inf_{
\begin{array}{c}
\tilde \lambda_i \geq 0, i\in \I(\bar w),\\
\sum_{i\in \I(\bar w)}\tilde \lambda_i = 1,\\
\tilde \mu_j \geq 0, j \in \J(\bar w),\\
\tilde \mu_j = 0, j \in \J\setminus \J(\bar w)
\end{array}}
\left\|\sum_{i\in \I(\bar w)}\tilde \lambda_i c_i +E\transp \tilde
v - \tilde\mu\right\|}.
$$
Recalling the notation of Section~1 allows us to reduce the latter
expression to the one in the lemma formulation and thus finish the
proof. $\h$\vspace*{0.05in}

Combining the obtained lemmas with the definitions above, we can
now complete the alternative proof of the condition measure
formula in Theorem~\ref{th3}.\\[1ex]
{\bf Proof of Theorem~\ref{th3}.} For every $w \in
{\Omega}\setminus S$ choose the parameter $z_{w}$ as in
Lemma~\ref{lem:GlobIsLoc}, i.e.,  put $z_{w} = F( w)- \e_{w}$,
where $\e_{w}$ is taken from Lemma~\ref{lem:OptCoincide}. Then we
have
\begin{eqnarray*}
k_F &   =   &   \inf\big\{ k\ge 0\big|\; \dist (w; S)\leq k F(w)
 \quad
                \mbox{for all }\; w\in {\Omega}\setminus S\big\}\qquad
                \big(\mbox{by definition \eqref{eq:SimplK}}\big)\\
    &   =   &   \sup_{w\in {\Omega}\setminus S} \frac{\dist(w; S)}{F(w)}\qquad \mbox
    {(as ${\Omega}\setminus S
    \ne \emptyset$)}\\
    &   =   &   \sup_{w\in {\Omega}\setminus S} \frac{\dist\big(w; S(z_w)\big)}{F(w) - z_w}\qquad
                \mbox{(by Lemma~\ref{lem:GlobIsLoc})}\\
    &   =   &   \sup_{w\in {\Omega}\setminus S} \frac{\Tilde V_{z_w}(w)}{F(w) - z_w}\qquad
                \big(\mbox{by the definition of ($P'_z$)\big)}\\
    &   =   &   \sup_{w\in {\Omega}\setminus S} \frac{V_{z_w}(w)}{F(w) - z_w}\qquad
                \mbox{(by Lemma~\ref{lem:OptCoincide} and the choice of $z_w$)}\\
    &   =   &   \sup_{w\in {\Omega}\setminus S}
                \frac{1}{\dist\big(0;\co_{i\in\I(w)} \{c_\ell\}+(\ker E)^\perp-
                \cone\big[\co_{j\in \J(w)} \{e_j\}\big]\big)} \qquad
                \mbox{(by Lemma~\ref{lem:slnSubprob})}.
\end{eqnarray*}
Letting $w=(x,y) \in \Omega$, observe finally that
\begin{eqnarray*}
\big\{c_\ell\big|\; \ell \in \I(w)\big\} = \big\{(a_i,b_k)\big|\;
i\in I(x),\; k\in K(y)\},\quad (\ker E)^\perp = \span \{{\bf 1}_n\}\times
\span \{{\bf 1}_m\},
\end{eqnarray*}
and $\k(A)= k_F$, which complete the proof of the theorem. $\h$\\

{\bf Acknowledgments.} The authors are indebted to two anonymous
referees for their valuable suggestions and remarks that allowed
us to essentially improve the original presentation. We also
gratefully acknowledge helpful discussions with Arkadi Nemirovski
and Yurii Nesterov on the topics and results of this paper.

\end{document}